# Improved Uncooperative Spacecraft Maneuver Detection with Space-based Optical Observations

Xuejian MAO[a]; Pei LIU[a]; Pei CHEN[a,*]

[a] *School of Astronautics, Beihang University, 102206 Beijing, China*

**Abstract**

Building and maintaining a space object catalog is necessary for space situational awareness. To realize this, one great challenge is uncooperative spacecraft maneuver detection because unknown maneuver events can lead to deviated orbital predictions and losses of tracking. Nowadays, more and more spacecraft equip electric propulsion and perform long-duration maneuvers to realize orbital transfer. Previous studies have investigated impulsive maneuver detection with space surveillance data. But, the developed method does not suffice for cases where maneuver durations are long. In this study, an improved uncooperative spacecraft maneuver detection method with space-based optical observations is proposed. Instead of a sudden maneuver event, the maneuver duration is considered. The maneuver starting/ending times are estimated along with the thrust acceleration vector. The angular residuals of nonlinear least square estimates are used to judge whether a maneuver policy could be a potential solution. The global optimum maneuver policy is chosen from multiple local minima according to the minimum-fuel principle. It is demonstrated that the maneuver duration is poorly observable if the thrust is along the orbital normal direction attributed to the nature of orbital dynamics. Real maneuver data of the Senitnel-3A spacecraft and the Senitnel-6A spacecraft is used to test the capability of the developed method.

## 1. Introduction

The number of space objects has been increasing explosively in recent years [1]. About 28000 space objects larger than 10 cm have been tracked by the U.S. space surveillance network till October 2024 [2]. Especially, the active

---

[*] Authors' addresses: Xuejian Mao, Pei Liu, and Pei Chen are with the School of Astronautics, Beihang University, 102206 Beijing, China, E-mail: (maoxuejian@buaa.edu.cn; zy2315116@buaa.edu.cn; chenpei@buaa.edu.cn). Tel.: +86 10 82316535. (Corresponding author: Pei Chen)

satellites have doubled in the last five years attributed to the development of mega constellations. The congestion in space greatly increases the risk of collision between space objects [3]. Frequent maneuvers are performed to avoid potential collisions. SpaceX reported that the Starlink spacecraft performed 49384 collision avoidance maneuvers during the period between December 1, 2023 and May 31, 2024 [4]. In the meantime, large amount of maneuvers makes the space traffic difficult to track and prediction, and brings more threats to operational spacecraft [5]. To make space environment safe and sustainable, maneuver detection and characterization of uncooperative spacecraft is critical in the current and upcoming space eras [6−8]. Tracking space objects is mainly through ground-based radar observations, ground-based optical observations, and space-based optical observations. In recent years, space-based space surveillance (SBSS) has gathered worldwide interests because space-based optical observations have great advantages in terms of temporal and spatial coverage compared to the ground-based observations [9−11]. For timely maneuver detection with space-based optical observations, one great difficulty lies in that the orbit can hardly be determined from single optical tracklet (a tracklet refers to a collection of consecutive observations in time), because optical observations only provide angular information and no range information [7,12].

Early studies mostly determine the maneuver information by comparing the post-maneuver orbit to the pre-maneuver orbit. The pre-maneuver orbit is determined from historical observations. The post-maneuver orbit can be individually solved with post-maneuver observations, where the initial orbit determination followed by the batch least square (BLS) procedure are generally utilized [13], or can be derived from the pre-maneuver orbit through sequential filters [5,14]. In recent years, various sequential filtering techniques for post-maneuver orbit determination have been developed. Kelecy et al. [15] used the extended Kalman filter (EKF) within the Orbit Determination Tool Kit (ODTK) for low-thrust maneuver detection. Goff et al. [16,17] applied variable dimension filters for orbit estimation of a continuously thrusting spacecraft and investigated the adaptive EKF and adaptive unscented Kalman filter (UKF) for noncooperative satellite maneuver reconstruction. Co et al. [18] represented the maneuver acceleration history by the sum of trigonometric series. The thrust Fourier coefficients and the orbital state are jointly estimated. Guang et al. [19] designed a variable structure estimator where the maneuver information is detected by a maneuver observer and fed into the nominal Kalman filter as compensation. Montilla et al. [20] developed an adaptive UKF for maneuver detection. The state covariance is adaptively changed so that the uncertainty of the state can include the effects of maneuver. In brief, sequential filter based methods address maneuvers as uncertain dynamics in the state propagation process, thus does not have a strong physical significance for arbitrary systems [6].

In realistic scenarios, the spacecraft maneuver control policy is generally designed to be optimal or suboptimal in terms of fuel or time consumption [21]. Motivated by this fact, Holzinger et al. [22] proposed to use the control usage as a metric to correlate object observations, detect maneuvers, and characterize maneuvers. Singh et al. [23] used the fuel cost (i.e., the total velocity increment) as the control distance metric and used nonlinear programming to solve the resulting optimal control problem. Lubey et al. [24] furtherly modified the control distance metric to incorporate the boundary state uncertainty. Compared to the sequential filter based maneuver detection methods, these optimal control based maneuver detection methods require less post-maneuver observations because of their rigorous physical meanings. Nonetheless, because the maneuver capabilities of uncooperative spacecraft are not preliminarily known, the estimated maneuver history may be significantly different from the true maneuver history. For example, a fuel-optimal orbital transfer comprises several impulsive maneuvers if the magnitude of thrust forces is not restricted. However, in cases where spacecraft equip low thrust propulsions, the fuel-optimal orbital transfer is the bang-bang control because the magnitude of thrust forces is limited [21]. Pastor et al. [7] developed a simplified maneuver detection methodology focusing on impulsive maneuvers. The impulsive maneuvers, i.e., sudden velocity changes, are estimated by searching for the local minimum of the weighted root mean square (WRMS) of observation residuals. One local minimum with the least control effort is selected as the final maneuver estimation. Porcelli et al. [8] employed a similar scheme to Pastor et al. [7] for LEO satellite maneuver detection with radar measurements. A novel dynamic model for maneuver detection was furtherly developed to address the strong perturbations of LEO satellites. Tests with real data show that the impulsive maneuver model cannot address long-duration maneuvers and can lead to biased maneuver estimations and tracklet correlation failures.

The main objective of maneuver detections is to answer the following two questions: 1) whether the spacecraft maneuvered, and 2) how the spacecraft maneuvered? Exiting methods can generally tackle the first question well. However, different methods approximate the maneuver history by different models. The maneuver estimation results with different methods remain diverse [25]. The method developed by Porcelli et al. [8] and Pastor et al. [7] focuses on impulsive maneuver detection and has been validated with real data, but cannot suffice for long-duration maneuver detection. Nowadays, more and more spacecraft employ electric propulsions and orbital maneuvers would last for long durations [26,27]. The present study considers effects of maneuver duration and develops a long-duration maneuver detection method with space-based optical observation. The starting time and ending time of the maneuver are used as optimization variables. With different combinations of maneuver starting time and maneuver ending time, an ensemble of post-maneuver orbits is computed with nonlinear least square estimator. Post-maneuver orbits with acceptable angular residuals are selected out as candidates. The candidate with minimum fuel cost is chosen as the final estimate. The real flight data of the Sentinel-3A satellite and the Sentinel-6A spacecraft are used to test the proposed maneuver detection method.

The reminder of this paper is organized as follows. Firstly, the space-based space surveillance system is briefly introduced. The mathematical model of space-based optical observation is built. The orbital evolution of a spacecraft implementing maneuvers is described. Secondly, the developed maneuver detection method is presented. The observability of maneuver parameters is briefly analyzed with aid of relative orbital motion. Thirdly, real flight data of the Sentinel-3A spacecraft and the Sentinel-6A spacecraft are used for tests. The capability and performance of the developed methodology are discussed. The conclusion of this study is drawn in the final section.

## 2. System Description

*2.1. Space-based Optical Observation*

The space-based optical observations are provided by space-based visible sensors onboard spacecraft platforms. The space-based visible sensor along with its spacecraft platform constitutes a space-based space surveillance (SBSS) system. Representative examples include the U.S. Air Force's SBV sensor on the Midcourse Space Experiment (MSX) spacecraft, the SBSS Block 10 System built by Boeing and Ball Aerospace, and the Canadian Sapphire and NEOSSat satellites [28]. These SBSS systems have greatly increased the revisit rates on space objects without the limitations inherent in ground systems.

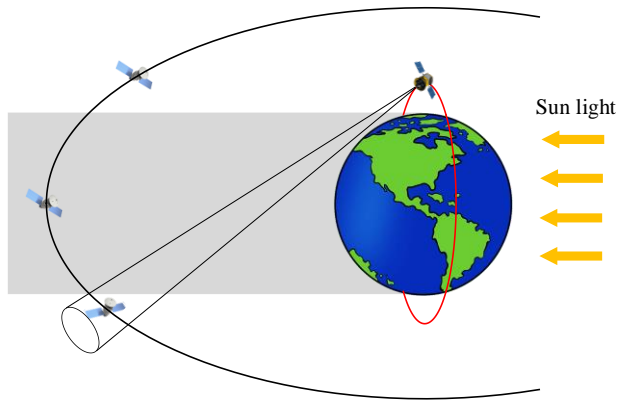

**Fig. 1.** SBSS scenario performing tasked tracking in the GEO belt

In a SBSS scenario, the SBV and its satellite platform are usually placed at a LEO, more generally, a Sun-synchronous orbit. The SBV is servo-controlled to point to the anti-Sun direction and stares at the target region. Figure 1 is a typical SBSS scenario performing tasked tracking in the GEO belt. The SBV receives the photons reflected off of a space object when this space object is illuminated by the Sun light. Based on the astrometric reduction against a background map of stars, the triangulation of the observed object with respect to catalogued stars can be performed, which derives the right ascension (RA) and declination (DEC) measurements in the Earth-centered inertial (ECI) coordinate frame, namely the space-based optical observations [29]. The RA and DEC measurement models are as follows

$$\alpha = \text{atan2}(y_t - y_s, x_t - x_s) + w_\alpha$$
$$\delta = \text{asin}\left(\frac{z_t - z_s}{\rho}\right) + w_\delta \quad (1)$$

with

$$\rho = \|\mathbf{r}_t - \mathbf{r}_s\| \quad (2)$$

where $\alpha$ and $\delta$ denote the RA and DEC respectively, atan2(,) is the usual two-argument arctangent function, $\|\cdot\|$ denotes the Euclidean norm, $\mathbf{r}_s = [x_s, y_s, z_s]^T$ is inertial position of the SBV, $\mathbf{r}_t = [x_t, y_t, z_t]^T$ is the inertial position of a space object, and $w_\alpha$ and $w_\delta$ are the measurement noise of the RA and DEC respectively. The superscript $T$ denotes the matrix transpose operation. In practical cases the inertial position of the SBV can be provided by on-ground precise orbit determination using GNSS measurements and has centimeter level accuracy, thus is assumed known in this study.

Equation (1) is a simplified model without considering dynamic effects including light time and aberration [30]. The measurement deviation caused by light time is due to the finite velocity of light, i.e., the time at which photons are reflected off of a space object differs from that at which they are received by the SBV. The light time is calculated by

$$\lambda \approx \frac{\rho}{c} \quad (3)$$

where $c = 3 \times 10^8$ m/s is the light velocity. The aberration is caused by the velocity of the SBV, i.e., there is a change in the direction of the light coming from a space object due to the motion of the LEO observer. With considering both light time and aberration effects, the line-of-sight vector of the GEO object with respect to the SBV is

$$\vec{n} = \frac{\mathbf{r}_t - \mathbf{r}_s - \frac{\rho}{c}\dot{\mathbf{r}}_t + \frac{\tilde{\rho}}{c}\dot{\mathbf{r}}_s}{\left\|\mathbf{r}_t - \mathbf{r}_s - \frac{\rho}{c}\dot{\mathbf{r}}_t + \frac{\tilde{\rho}}{c}\dot{\mathbf{r}}_s\right\|} \tag{4}$$

with

$$\tilde{\rho} = \left\|\mathbf{r}_t - \mathbf{r}_s - \frac{\rho}{c}\dot{\mathbf{r}}_t\right\| \tag{5}$$

Because $\tilde{\rho} \approx \rho$, the equivalent angle measurements of Eq. (4) is

$$\alpha \approx \mathrm{atan2}\left(y_t - y_s + \frac{\rho}{c}(\dot{y}_s - \dot{y}_t),\ x_t - x_s + \frac{\rho}{c}(\dot{x}_s - \dot{x}_t)\right) + w_\alpha$$

$$\delta \approx \mathrm{asin}\left(\frac{z_t - z_s + \frac{\rho}{c}(\dot{z}_s - \dot{z}_t)}{\breve{\rho}}\right) + w_\delta \tag{6}$$

with

$$\breve{\rho} = \left\|\mathbf{r}_t - \mathbf{r}_s - \frac{\rho}{c}\dot{\mathbf{r}}_t + \frac{\rho}{c}\dot{\mathbf{r}}_s\right\| \tag{7}$$

The partial derivatives of angular measurements with respect to the target's position and velocities are

$$\mathbf{H}_r = \frac{\partial[\alpha,\delta]^T}{\partial \mathbf{r}_t} \approx \begin{bmatrix} -\frac{y_t - y_s}{(x_t - x_s)^2 + (y_t - y_s)^2} & -\frac{(x_t - x_s)\cdot(z_t - z_s)}{\rho^2 \sqrt{(x_t - x_s)^2 + (y_t - y_s)^2}} \\ \frac{x_t - x_s}{(x_t - x_s)^2 + (y_t - y_s)^2} & -\frac{(y_t - y_s)\cdot(z_t - z_s)}{\rho^2 \sqrt{(x_t - x_s)^2 + (y_t - y_s)^2}} \\ 0 & \frac{\sqrt{(x_t - x_s)^2 + (y_t - y_s)^2}}{\rho^2} \end{bmatrix}^T \tag{8}$$

$$\mathbf{H}_v = \frac{\partial[\alpha,\delta]^T}{\partial \mathbf{v}_t} \approx \frac{\rho}{c}\mathbf{H}_r \tag{9}$$

Above partial derivative matrices will be used in the maneuver detection algorithm.

## 2.2. Spacecraft Maneuver Models

The orbital motion of a controlled spacecraft is governed by gravitational forces, perturbation forces, and control forces, which can be represented by the following differential equation expressed in the ECI frame:

$$\frac{\mathrm{d}}{\mathrm{d}t}\begin{bmatrix} \mathbf{r} \\ \mathbf{v} \end{bmatrix} = \begin{bmatrix} \mathbf{v} \\ \mathbf{g}(t,\mathbf{r}) + \mathbf{a}_p(t,\mathbf{r},\mathbf{v}) + \mathbf{f}(t) \end{bmatrix} \tag{10}$$

where $r$ and $v$ are the spacecraft's position and velocity respectively, $t$ denotes time. $g$ is the gravitational acceleration, $a_p$ is the perturbation forces including third-body attractions from the Moon and the Sun, atmospheric drag, solar radiation, etc. $f$ is the thrust acceleration.

If the magnitude of the thrust acceleration is large and the maneuver duration is very short, the maneuver can be regarded as impulsive, i.e., a sudden velocity change takes places at the maneuver epoch as follows:

$$\begin{aligned} r(t_M^+) &= r(t_M^-) \\ v(t_M^+) &= v(t_M^-) + \Delta v \end{aligned} \quad (11)$$

where $t_M$ is the maneuver epoch, and $\Delta v$ is the sudden velocity change. Equation (11) is a mathematical approximation of short-burn maneuver cases and has been proved valid for short-burn maneuver estimation [7,17].

When a large $\Delta v$ is required, the thruster need be fired for a long time. In addition, electric propulsions are widely used in nowadays. The thrust forces with electric propulsions are usually small. As a result, the maneuver duration is long and Eq. (11) cannot be applied. In practical long-burn maneuver missions, the thrust forces mostly stay constant in the orbital frame. This approach is only suboptimal but is easy to implement, thus is widely adopted. For example, in orbital maintenance maneuvers thrust forces generally stay in the direction of spacecraft motion [21]. The long-duration constant-thrust maneuver can be modelled as:

$$\begin{cases} \frac{d}{dt}\begin{bmatrix} r \\ v \end{bmatrix} = \begin{bmatrix} v \\ g(t,r) + a_p(t,r,v) \end{bmatrix}, & t \leq t_f \text{ or } t > t_b \\ \frac{d}{dt}\begin{bmatrix} r \\ v \end{bmatrix} = \begin{bmatrix} v \\ g(t,r) + a_p(t,r,v) + C(r,v) \cdot u \end{bmatrix}, & t_b < t \leq t_f \end{cases} \quad (12)$$

where $t_b$ and $t_f$ are the starting and ending times of the maneuver, respectively, $u$ is the constant thrust acceleration in the orbital frame, and $C(r,v)$ is the coordinate transformation matrix from the orbital frame to the ECI frame. The orbital frame, also known as the Vehicle Velocity, Local Horizontal (VVLH) frame, is defined as follows: the Z-axis is along the negative position vector, the Y-axis is along the negative orbit normal, and the X-axis is toward velocity.

For the problem of maneuver detection, the maneuver type of an uncooperative spacecraft is not preliminarily known. In cases where maneuver durations are not quite long, true maneuver histories can be approximated by sudden velocity changes. Herein, the discrepancy between long-duration maneuvers and corresponding impulsive maneuvers is quantitatively assessed. A circular low Earth orbit with 500-km height is simulated under the two-body dynamics. The orbital period is about 95 min. The orbital states of a long-duration maneuvering spacecraft are obtained by integrating the ordinary differential equation (12). Another orbit with the same initial orbital sates and a corresponding impulsive maneuver is also generated. The impulsive maneuver takes place at the middle epoch of the long-duration maneuver and has the magnitude of $(t_b - t_f)\|\bm{u}\|$ (denoted by $\Delta V$).

The 3D position differences between two set of orbital states are shown in Fig. 2, where $\Delta t = t_b - t_f$ is the maneuver duration. The magnitude of thrust acceleration is set to a typical value of $10^{-3}$ m/s². Three representative cases are investigated including constant thrust along the in-track direction, constant thrust along the radial direction, and constant thrust along the orbital normal direction. It is seen that the 3D position differences between two set of orbital states are below 1 m after 24-h propagation if the maneuver duration is as small as 60 s. For the case with the constant thrust along the in-track direction, with maneuver durations of 300 s, 600 s, 1200 s, 1800 s, and 2400 s, the average 3D position differences between two set of orbital states during 1-day propagation are 3.9 m, 30.6 m, 241.2 m, 791.8 m, and 2506.2 m, respectively. For the case with the constant thrust along the radial direction, with maneuver durations of 300 s, 600 s, 1200 s, 1800 s, and 2400 s, the average 3D position differences between two set of orbital states during 1-day propagation are 1.9 m, 15.2 m, 119.8 m, 393.2 m, and 1520.0 m, respectively. For the case with the constant thrust along the orbital normal direction, with maneuver durations of 300 s, 600 s, 1200 s, 1800 s, and 2400 s, the average 3D position differences between two set of orbital states during 1-day propagation are 1.5 m, 7.6 m, 52.9 m, 168.9 m, and 526.5 m, respectively. In addition, the 3D position differences are almost proportional with respect to the magnitude of thrust acceleration. For example, if the magnitude of thrust acceleration is increased to $10^{-2}$ m/s², the average 3D position difference during 1-day propagation is 306.9 m with 600-s maneuver duration for the case with the constant thrust along the in-track direction. The present study mainly focuses on normal orbital adjustment maneuvers. Large orbital transfers, such as LEO-to-GEO transfer, are not investigated.

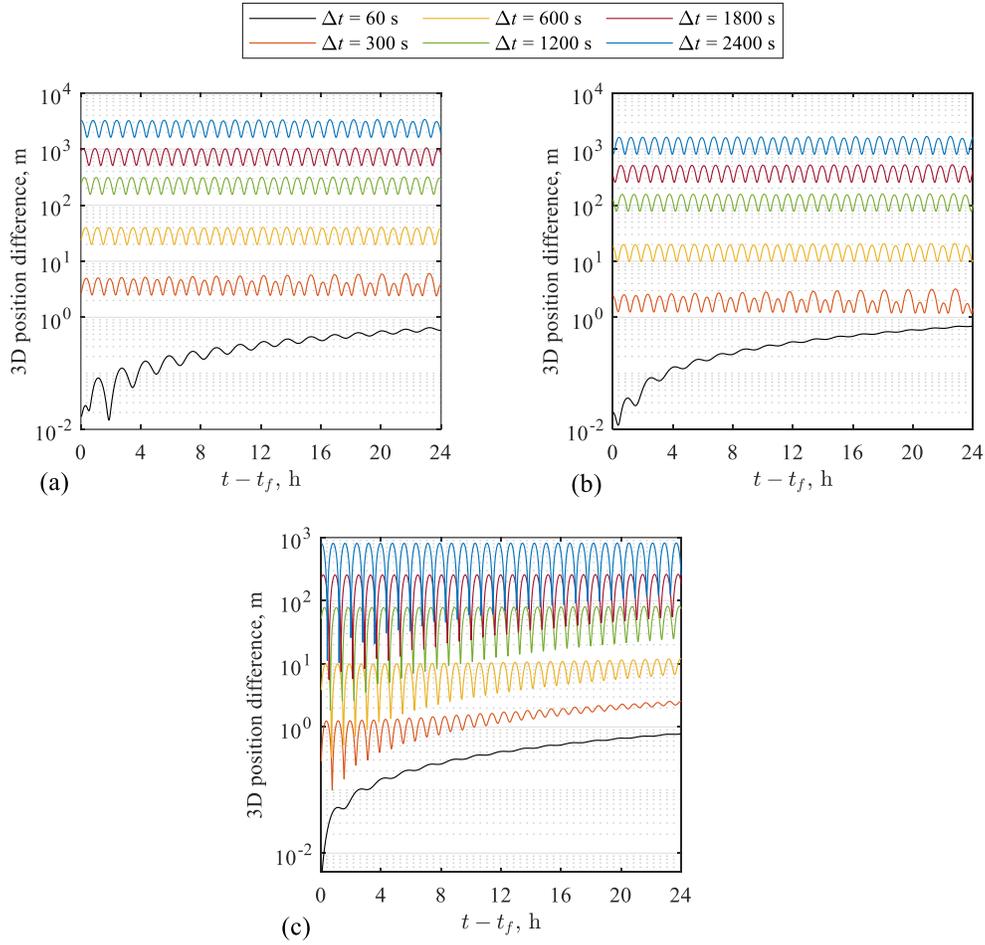

**Fig. 2.** 3D position difference between two orbits undergoing long-duration constant-thrust maneuver and corresponding impulsive maneuver. (a) constant thrust along the in-track direction with $\boldsymbol{u} = \begin{bmatrix} 10^{-3}, 0, 0 \end{bmatrix}^T$ m/s$^2$. (b) constant thrust along the radial direction with $\boldsymbol{u} = \begin{bmatrix} 0, 0, 10^{-3} \end{bmatrix}^T$ m/s$^2$. (c) constant thrust along the orbital normal direction with $\boldsymbol{u} = \begin{bmatrix} 0, 10^{-3}, 0 \end{bmatrix}^T$ m/s$^2$.

## 3. Maneuver Detection Methodology

### *3.1. Post-maneuver Tracklet Correlation*

For maneuver detection with space-based optical observations, one prior work is to correlate post-maneuver tracklets to pre-maneuver orbits, i.e., the tracklet-to-orbit (T2O) problem. Several methods have been proposed for the T2O problem and proven to be effective [8, 22, 31]. However, although T2O methods may answer the question whether a spacecraft maneuvered, the maneuver estimation results with single post-maneuver tracklet are not reliable [7,8]. If two or more post-maneuver tracklets are used to detect maneuvers, these tracklets must be correlated to each other, i.e., the tracklet-to-tracklet (T2T) problem. Representative methods for the T2T problem include the IVP method and the BVP method [32,33]. In this study, two or more post-maneuver tracklets are correlated to each other at first in the working flow of maneuver detection. Then, an initial estimate of the post-maneuver orbit is obtained from these post-maneuver tracklets with a least square estimator. The correlation between the pre-maneuver orbit and post-maneuver tracklets is achieved through the correlation between the pre-maneuver orbit and the estimated post-maneuver orbit, i.e., the orbit-to-orbit (O2O) problem.

Denote the propagated position history of a pre-maneuver orbit by $\{r_1(t_1), r_1(t_2), \ldots, r_1(t_n)\}$. The covariance matrices of positions are $\{P_1(t_1), P_1(t_2), \ldots, P_1(t_n)\}$. Denote the propagated position history of an estimated post-maneuver orbit by $\{r_2(t_1), r_2(t_2), \ldots, r_2(t_n)\}$. The covariance matrices of positions are $\{P_2(t_1), P_2(t_2), \ldots, P_2(t_n)\}$. In impulsive maneuver cases, the O2O problem can be simply solved by determining whether two orbits can intersect at certain epochs as follows [34]:

$$\exists i \in \{1, 2, \ldots, n\}, \quad \text{s.t.} \quad \chi(t_i) \leq \chi_{\max}$$
$$\text{with} \quad \chi(t_i) = \sqrt{(r_2(t_i) - r_1(t_i))^T (P_2(t_i) + P_1(t_i))^{-1} (r_2(t_i) - r_1(t_i))} \tag{13}$$

where $\chi(t_i)$ is the Mahalanobis distance between two positions at epoch $t_i$, $\chi_{\max}$ is the threshold of Mahalanobis distance and can be set to 3.38, which encompasses approximately 99% of the Mahalanobis distance cumulative distribution function within the three-dimensional space. In applications, the Mahalanobis distances at all sampling epochs are calculated to check whether the intersection condition is satisfied.

In long-duration thrust cases, the post-maneuver orbit doesn't physically intersect with the pre-maneuver orbit. According to results in Fig. 2, the post-maneuver orbit is around a virtual orbit that can be transferred into from the pre-maneuver orbit with an impulsive maneuver at the middle epoch of the maneuver duration. Thus, Eq. (13) can be adapted into

$$\exists i \in \{1, 2, \ldots, n\} \text{ and } \Delta r \in \mathbb{R}^3, \quad \text{s.t.} \quad \begin{cases} \chi(t_i) \leq \chi_{\max}^2 \\ \Delta r^T \Delta r \leq d^2 \end{cases}$$
$$\text{with} \quad \chi(t_i) = \sqrt{(r_2(t_i) - r_1(t_i) - \Delta r)^T (P_2(t_i) + P_1(t_i))^{-1} (r_2(t_i) - r_1(t_i) - \Delta r)} \tag{14}$$

where $\Delta \boldsymbol{r}$ represents the position difference between the post-maneuver orbit and the virtual orbit, and $d$ is the threshold representing the maximum of the 3D position difference and is set to 10 km in this study. The solution procedure for Eq. (14) is somewhat sophisticated and is provided in the Appendix.

A numerical example is developed to illustrate the differences between Eq. (13) and Eq. (14). The spacecraft undergoes a constant thrust along the in-track direction with $\boldsymbol{u} = \begin{bmatrix} 10^{-3}, 0, 0 \end{bmatrix}^T$ m/s². The maneuver duration is 1800 s. Two post-maneuver tracklets are used to obtain an initial estimate of the post-maneuver orbit. The computed Mahalanobis distances are shown in Fig. 3. It is seen that the smallest value of Mahalanobis distance with Eq. (13) is 4.987 and is larger than the threshold. Whereas, the computed Mahalanobis distances with Eq. (14) are below the threshold when their epochs are near the maneuver time. Therefore, Eq. (14) is capable for O2O correlation in long-duration maneuver cases. In addition, the Mahalanobis distances indicate the range of the middle epoch of the maneuver duration.

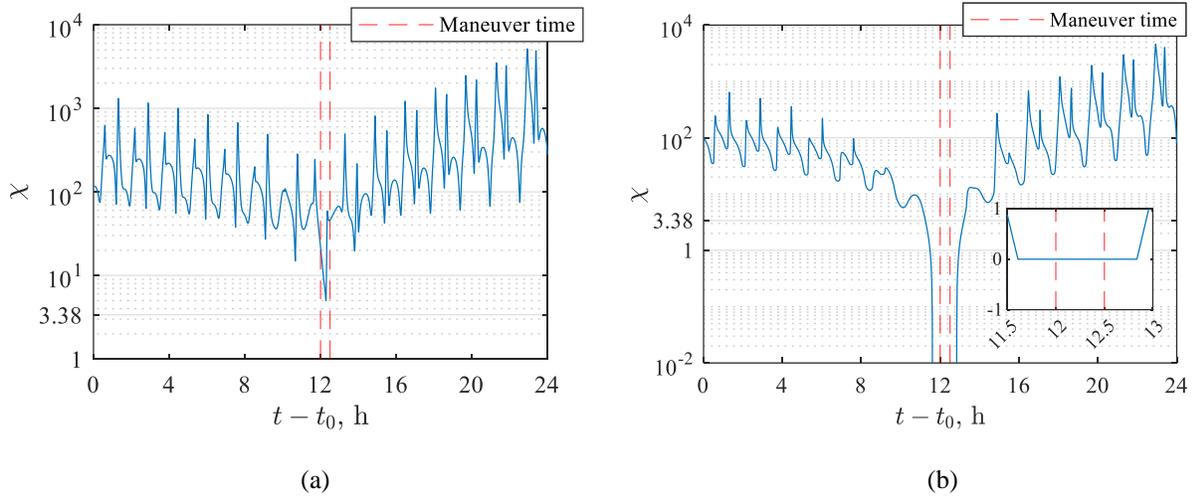

**Fig. 3.** An example of Mahalanobis distances with different methods. (a) computed with Eq. (13). (b) computed with Eq. (14)

*3.2. Long-duration Constant-thrust Maneuver Detection*

After post-maneuver observations are correlated to the pre-maneuver orbit, the estimate of the maneuver starting/ending time and thrust acceleration is triggered. All post-maneuver optical observations are stacked into a measurement vector as

$$Y = \begin{bmatrix} y(t_1) \\ y(t_2) \\ \vdots \\ y(t_n) \end{bmatrix} \tag{15}$$

where $y(t_k) = [\alpha(t_k), \delta(t_k)]^T$ is the angular measurement at epoch $t_k$ as described in Eqs. (1) and (6). Give a pair of maneuver starting/ending time guesses $\tilde{t}_b$ and $\tilde{t}_f$, the optimal thrust acceleration under this pair of maneuver time guesses is computed with the nonlinear least square estimator in a iterative manner as follows [35]

$$u \Leftarrow u + \left[ \left( \frac{\partial Y}{\partial u^T} \right)^T R_Y^{-1} \left( \frac{\partial Y}{\partial u^T} \right) \right]^{-1} \left( \frac{\partial Y}{\partial u^T} \right)^T R_Y^{-1} \left( Y - Y^-(\tilde{t}_b, \tilde{t}_f, u) \right) \tag{16}$$

where $\frac{\partial Y}{\partial u^T}$ is the partial derivative matrix of the stacked measurement vector with respect to the thrust acceleration, $R_Y$ is the covariance matrix of $Y$, and $Y^-(\tilde{t}_b, \tilde{t}_f, u)$ is the predicted measurements with mathematical models in Eq. (6). $\frac{\partial Y}{\partial u^T}$ is calculated with chain rule as

$$\frac{\partial Y}{\partial u^T} = \begin{bmatrix} \frac{\partial y(t_1)}{\partial u^T} \\ \frac{\partial y(t_2)}{\partial u^T} \\ \vdots \\ \frac{\partial y(t_n)}{\partial u^T} \end{bmatrix} \tag{17}$$

with

$$\frac{\partial y(t_k)}{\partial u^T} = [H_r(t_k) \quad H_v(t_k)] \Phi(t_k, \tilde{t}_f) S(\tilde{t}_f, \tilde{t}_b) \tag{18}$$

where $\Phi(t_k, \tilde{t}_f)$ is the orbital state transition matrix from $\tilde{t}_f$ to $t_k$. $S(\tilde{t}_f, \tilde{t}_b)$ is the sensitivity matrix of the orbital state at $\tilde{t}_f$ with respect to the thrust acceleration $u$. Denote the initial epoch of orbital propagation by $t_0$. $\Phi(t_k, \tilde{t}_f)$ is calculated by

$$\Phi(t_k, \tilde{t}_f) = \Phi(t_k, t_0) \Phi^{-1}(\tilde{t}_f, t_0) \tag{19}$$

$$S(\tilde{t}_f, \tilde{t}_b) = S(\tilde{t}_f, t_0) - \Phi(\tilde{t}_f, \tilde{t}_b) S(\tilde{t}_b, t_0) \tag{20}$$

$\boldsymbol{\Phi}(t,t_0)$ and $\boldsymbol{S}(t,t_0)$ are simultaneously obtained by integrating the following ordinary differential equations [31]:

$$\dot{\boldsymbol{\Phi}} = \boldsymbol{A}\boldsymbol{\Phi} \tag{21}$$

$$\dot{\boldsymbol{S}} = \boldsymbol{A}\boldsymbol{S} + \begin{bmatrix} \boldsymbol{0}_{3\times 3} \\ \boldsymbol{C}(\boldsymbol{r},\boldsymbol{v}) \end{bmatrix} \tag{22}$$

with

$$\boldsymbol{A} = \begin{bmatrix} \boldsymbol{0}_{3\times 3} & \boldsymbol{I}_{3\times 3} \\ \dfrac{\partial \boldsymbol{g}}{\partial \boldsymbol{r}^T} & \boldsymbol{0}_{3\times 3} \end{bmatrix} \tag{23}$$

$$\boldsymbol{\Phi}(t_0,t_0) = \boldsymbol{I}_{6\times 6} \; ; \; \boldsymbol{S}(t_0,t_0) = \boldsymbol{0}_{6\times 3} \tag{24}$$

When $\tilde{t}_f$ is very close to $\tilde{t}_b$, we have

$$\boldsymbol{S}(\tilde{t}_f,\tilde{t}_b) \approx \begin{bmatrix} \boldsymbol{0}_{3\times 3} \\ \boldsymbol{C}(\boldsymbol{r}(\tilde{t}_b),\boldsymbol{v}(\tilde{t}_b)) \end{bmatrix}(\tilde{t}_f - \tilde{t}_b) \tag{25}$$

Substitute Eq. (25) into Eqs. (18) and (16). It is seen that the algorithm estimates the sudden velocity change and becomes the impulsive maneuver detection method developed by Pastor et al. [7] and Porcelli et al. [8].

The performance index of above nonlinear least square solution is given as [31]

$$J(\tilde{t}_b,\tilde{t}_f) = \sqrt{\dfrac{1}{2n}\left(\boldsymbol{Y}-\boldsymbol{Y}^{-}(\tilde{t}_b,\tilde{t}_f,\hat{\boldsymbol{u}})\right)^T \boldsymbol{R}_Y^{-1}\left(\boldsymbol{Y}-\boldsymbol{Y}^{-}(\tilde{t}_b,\tilde{t}_f,\hat{\boldsymbol{u}})\right)} \tag{26}$$

It is seen that $J$ is the WRMS of angular residuals. For any pair of maneuver starting/ending time guesses, there is a corresponding performance index of the least square solution. Therefore, a pair of maneuver starting/ending time guess and the derived performance indexes constitute a one-way mapping. Each pair of maneuver starting/ending time guess and the derived thrust acceleration from Eq. (16) together represent a guess of the post-maneuver orbit. The true post-maneuver orbit naturally has good agreement with post-maneuver optical tracklets and the corresponding performance index should be close to 1 according to the law of large numbers. When the maneuver starting/ending time guess is equal to the actual maneuver starting/ending time, the resulting post-maneuver orbit guess will be very close to the true post-maneuver orbit from the statistical sense of the least square method. Therefore, a threshold of the performance index can be set beforehand. Only maneuver starting/ending time guesses having performance indexes smaller than this threshold are reserved as candidate solutions. Among all candidate solutions, the one with the minimum fuel cost (i.e., total velocity increment) is chosen as the final result according to the optimal control principle.

The scheme for uncooperative spacecraft maneuver detection is summarized in Fig. 3.

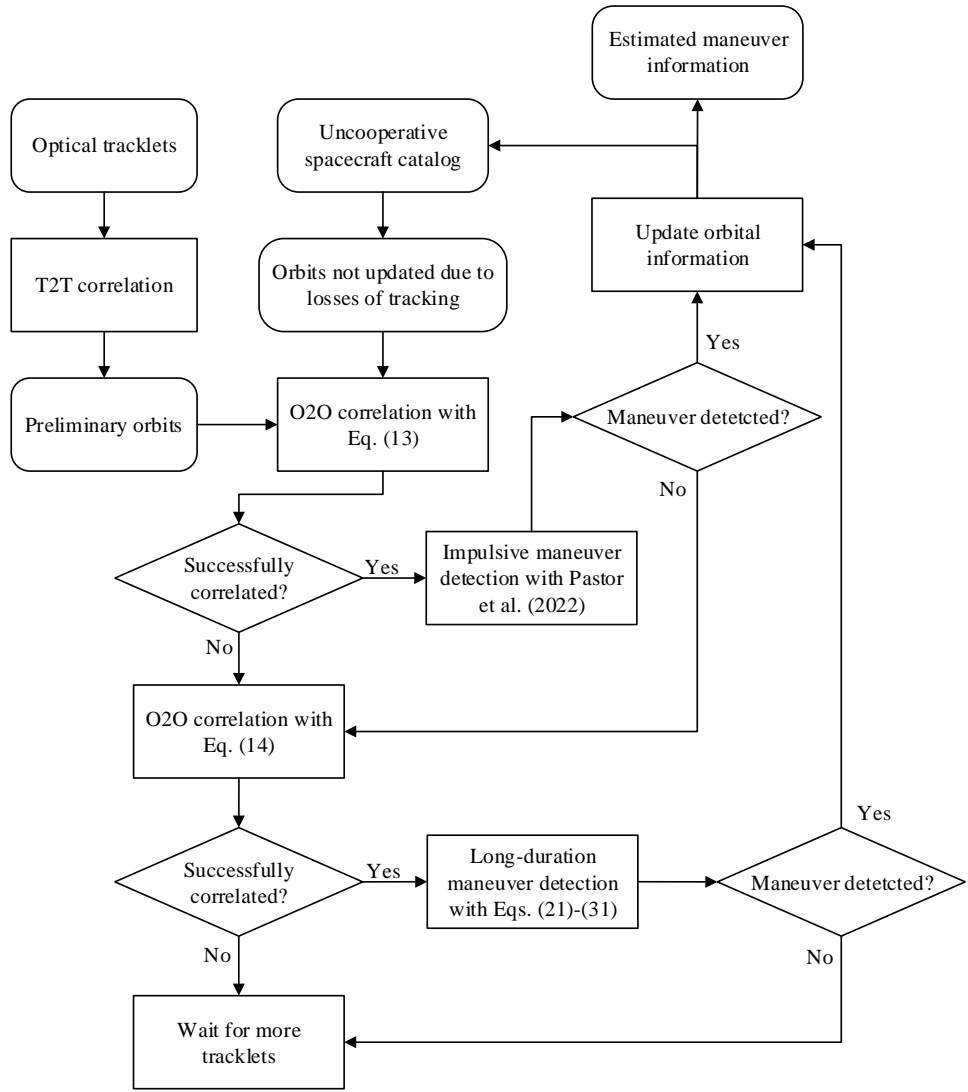

**Fig. 4.** The scheme for uncooperative spacecraft maneuver detection

*3.3. Observability Analysis with Virtual Position Measurements*

The system observability describes whether system states can be estimated from system output [36]. In the long-duration constant-thrust maneuver detection problem, there are five variables estimated with post-maneuver observations, including the maneuver starting time, the maneuver ending time, and the 3-dimensional thrust acceleration vector. When sufficient post-maneuver observations are collected, an accurate post-maneuver orbit can be estimated, which provides orbital information with six degrees of freedom. That is to say, long-duration constant-thrust maneuver detection can be roughly regarded as solving 6 nonlinear equations containing 5 unknown variables.

Intuitively speaking, a unique solution or no solution exists for this problem. But, there might be multiple solutions due to system nonlinearity, i.e., the system is unobservable.

Herein, we use the sequential position data of the post-maneuver orbit as virtual measurements instead of angular measurements for simplicity. With the developed maneuver detection algorithm, the thrust accelerations under different combinations of maneuver starting/ending time guesses are estimated. The position residuals suggest whether the pre-maneuver orbit can be transferred into the post-maneuver orbit with the estimated maneuver policy. The simulation conditions are as follows. Only the central force gravity field is considered. The pre-maneuver orbit is circular and has an orbital height of 500 km. The maneuver starting time is 12 h away from the initial epoch and the maneuver duration is 1200 s. The magnitude of thrust acceleration is set to $5\times10^{-3}$ m/s$^2$. Three different cases are investigated including constant thrust along the in-track direction, constant thrust along the radial direction, and constant thrust along the orbital normal direction. Sequential position data with elapsed time of 24 h ~ 48 h away from the initial epoch is used as virtual measurements and has a sampling interval of 60 s.

The root mean squares (RMS) of 3D position residuals is shown in Fig. 5. It is evidently seen that multiple local minima exist due to the periodicity of orbital motion. For all three cases, the global minima are obtained with the true maneuver starting/ending times. For the case with constant thrust along the in-track direction, local minima are generally larger than 1000 m except for the global minima. For the case with constant thrust along the radial direction, local minima are generally larger than 100 m except for the global minima. For the case with constant thrust along the orbital normal direction, there are numerous local minima smaller than 10 m, which indicates that the maneuver parameters are poorly observable in this case. In addition, it is shown that these local minima are almost located in lines with slope of −1. These results suggest that the pre-maneuver orbit need be quite accurate and the number of post-maneuver observations need be adequate, otherwise the true maneuver starting/ending can hardly be distinguished from multiple local minima only according to observation residuals.

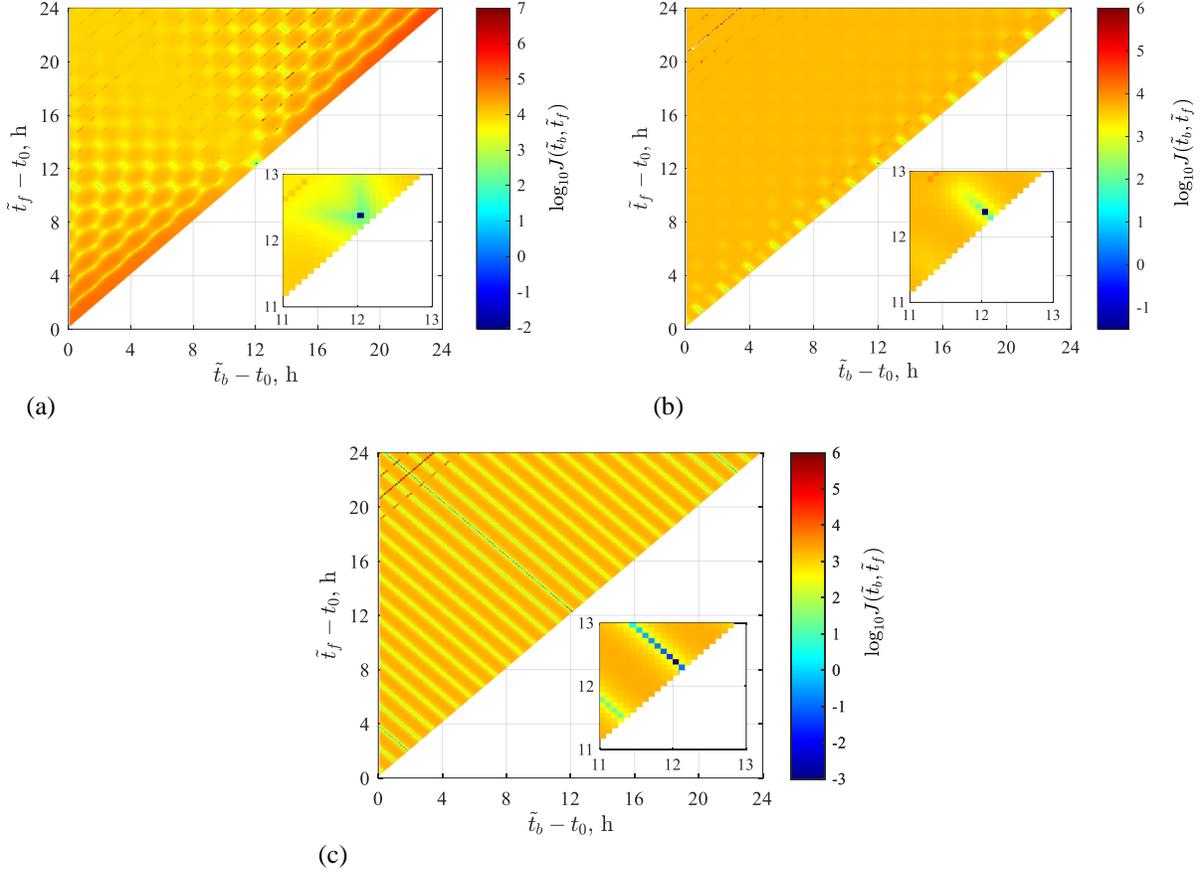

**Fig. 5.** RMS of 3D position residuals for maneuver detection with virtual position measurements. (a) constant thrust along the in-track direction. (b) constant thrust along the radial direction. (c) constant thrust along the orbital normal direction

The poor observability of the constant thrust along the orbital normal direction is unexpected and can be explained by relative orbital motion. A relative orbital model describes the relative motion of a deputy orbit with respect to the chief orbit. When the chief orbit is circular and the deputy orbit is near the chief orbit, the relative orbital motion is described by the Clohessy-Wiltshire (CW) equations [37]:

$$\begin{cases} \Delta \ddot{r}_x - 2n\Delta \dot{r}_z - 3n^2 \Delta r_x = u_x \\ \Delta \ddot{r}_y + n^2 \Delta r_y = u_y \\ \Delta \ddot{r}_z + 2n\Delta \dot{r}_x = u_z \end{cases} \qquad (27)$$

where $[\Delta r_x, \Delta r_y, \Delta r_z]^T$ is the relative position in the orbital frame, $n$ is the orbital rate of the chief orbit, and $[u_x, u_y, u_z]^T$ is the thrust acceleration of the deputy spacecraft. It is seen that the motion in the orbital normal direction is decoupled from those in the other directions. If there is no thrust acceleration, the motion in the orbital normal direction can be analytically expressed as

$$\begin{cases} \Delta r_y(t) = \dfrac{\Delta \dot{r}_y(t_0)}{n} \sin n(t-t_0) + \Delta r_y(t_0) \cos n(t-t_0) \\ \Delta \dot{r}_y(t) = \Delta \dot{r}_y(t_0) \cos n(t-t_0) - n \Delta r_y(t_0) \sin n(t-t_0) \end{cases} \quad (28)$$

With different initial conditions, Eq. (28) represents a cluster of trajectories in the $\Delta r_y$-$\Delta \dot{r}_y$ plane, as the solid ellipses in Fig. 5. If there is constant thrust acceleration in the orbital normal direction, Eq. (28) is adapted into

$$\begin{cases} \Delta r_y(t) = \dfrac{\Delta \dot{r}_y(t_0)}{n} \sin n(t-t_0) + \left[\Delta r_y(t_0) - \dfrac{u_y}{n^2}\right] \cos n(t-t_0) + \dfrac{u_y}{n^2} \\ \Delta \dot{r}_y(t) = \Delta \dot{r}_y(t_0) \cos n(t-t_0) - \left[n \Delta r_y(t_0) - \dfrac{u_y}{n}\right] \sin n(t-t_0) \end{cases} \quad (29)$$

In the case of spacecraft maneuver, the pre-maneuver orbit and the post-maneuver orbit are regarded as chief orbit and deputy orbit respectively, and the initial relative position and velocity are zeros. Equation (29) is changed into

$$\begin{cases} \Delta r_y(t) = -\dfrac{u_y}{n^2} \cos n(t-t_0) + \dfrac{u_y}{n^2} \\ \Delta \dot{r}_y(t) = \dfrac{u_y}{n} \sin n(t-t_0) \end{cases} \quad (30)$$

With different thrust acceleration, Eq. (30) also represents a cluster of trajectories in the $\Delta r_y$-$\Delta \dot{r}_y$ plane, as the dashed ellipses in Fig. 5. The poor observability of the constant thrust along the orbital normal direction is explained as follows. Under the true thrust acceleration, the spacecraft moves from point O to point A along the blue dashed arc. After the thruster is turned off at point A, the spacecraft is transferred into the target orbit and moves from point A to point B along the blue solid arc. However, the spacecraft can undergo another constant-thrust maneuver history and moves from point O to point B along the red dashed arc. If the thruster is turned off at point B, the spacecraft is transferred into the same target orbit. Therefore, multiple combinations of thrust acceleration and maneuver starting/ending times can achieve the same orbit transfer. Particularly, the spacecraft can be transferred into the target orbit with an impulsive maneuver, i.e., instantaneously moves from point O to point D.

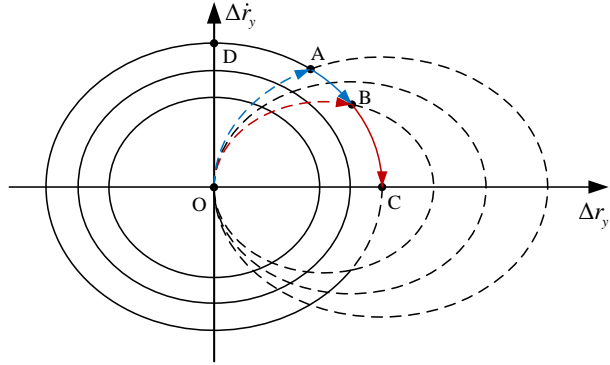

**Fig. 6.** Illustration of the poor observability of the constant thrust along the orbital normal direction

## 4. Results and Discussion

### *4.1. Data Preparation*

The developed maneuver detection algorithm is tested with the real flight data of the Sentinel-3A spacecraft and the Sentinel-6A spacecraft. The precise orbital states of these two spacecraft along with their maneuver histories can be downloaded from the International DORIS Service[†]. A spacecraft in a Sun-synchronous orbit is simulated as the SBSS platform and carries an optical sensor to provide space-based optical observations. The orbital states of the Sentinel-3A spacecraft, the Sentinel-6A spacecraft, and the simulated SBSS platform at the reference epoch 2020-12-13, 00:00:00.00 (UTC) are provided in Table 1. Two maneuver events are adopted for tests. The maneuver information is listed in Table 2. The information of the dynamical environment is as follows: The EGM2008 model truncated at the 20th degree and order for nonspherical gravitational perturbation, the NRLMSISE-00 model for the atmospheric drag, the DE430 ephemerides for the 3-body attractions from the Moon and the Sun, and the cannonball model for solar radiation pressure [37]. The optical sensor onboard the SBSS platform is assumed to have a wide FOV of 8°×8°. Angular observations are generated every 1 s with 5″ (1σ) uncertainties.

---

[†] https://cddis.nasa.gov/archive/doris/products/orbits/ssa

**Table 1.** Orbital states of the simulated SBSS platform, the Sentinel-3A spacecraft, and the Sentinel-6A spacecraft at 2020-12-13, 00:00:00.00 (UTC)

| Orbital parameters | SBSS platforms | Sentinel-3A | Sentinel-6A |
|---|---|---|---|
| Semi-major axis, km | 6888.580 | 7169.856 | 7706.232 |
| Eccentricity | $1 \times 10^{-4}$ | $5.918 \times 10^{-4}$ | $1.841 \times 10^{-3}$ |
| Inclination, ° | 97 | 98.715 | 66.037 |
| Right ascension of the ascending node, ° | 148 | 51.884 | 354.233 |
| Argument of the perigee, ° | 19 | 330.846 | 86.872 |
| Mean anomaly, ° | 275 | 132.233 | 296.094 |

**Table 2.** Maneuver information of the Sentinel-3A spacecraft and the Sentinel-6A spacecraft

| Maneuver parameters | Sentinel-3A | Sentinel-6A |
|---|---|---|
| Starting time (UTC) | 2020-12-16, 11:39:02 | 2020-12-14, 05:15:42 |
| Ending time (UTC) | 2020-12-16, 11:55:40 | 2020-12-14, 05:24:27 |
| Duration, s | 998 | 525 |
| Thrust acceleration vector (VVLH frame), mm/s² | $[1.318 \times 10^{-2}, -2.403, 4.121 \times 10^{-3}]$ | $[9.842, -9.647 \times 10^{-4}, 1.791 \times 10^{-2}]$ |
| $\Delta V$, m/s | 2.398 | 5.168 |

The initial epoch of the numerical scenario is set to 2020-12-13, 00:00:00.00 (UTC). During a 6-day simulation period, the obtained optical tracklets for the Sentinel-3A spacecraft and the Sentinel-6A spacecraft are shown in Fig. 7. For the Sentinel-3A spacecraft, a total of 15 optical tracklets are collected including 8 pre-maneuver tracklets and 7 post-maneuver tracklets. For the Sentinel-6A spacecraft, a total of 21 optical tracklets are collected including 7 pre-maneuver tracklets and 14 post-maneuver tracklets. The lengths of optical tracklets are 37 s on average. In a typical data processing flow of a SBSS system, pre-maneuver tracklets can be successfully correlated to orbits stored in the catalog, i.e., the pre-maneuver orbits [8]. The epoch of the last pre-maneuver observation is labeled as the latest update time of an orbit, denoted by $t_0$. After $t_0$, this orbit is not updated for a long time because of maneuver and the post-maneuver tracklets cannot be correlated to freely propagated pre-maneuver orbits. With the scheme in Fig. 3, a number of post-maneuver tracklets will be correlated to the pre-maneuver orbit under the hypothesis that the spacecraft performed maneuver. Then, the maneuver detection is triggered to confirm maneuver event and to estimate the maneuver information. Denote the earliest epoch of these post-maneuver tracklets by $t_1$. It is straightforward that the maneuver starting/ending times lie in the time span $[t_0,\ t_1]$. The developed maneuver detection method computes potential maneuver history with different combinations of maneuver starting/ending times within this time span. For simplicity, the pre-maneuver orbital states are assumed to be error-free as adequate pre-maneuver observations can usually be collected for precise orbit determination. In addition, the O2O correlation results indicate the range of the possible middle epochs of the maneuver duration. Another consideration is that the duration of single maneuver is generally no more than 1 h due to power limitation. These information reduces the number of pairs of maneuver starting/ending time guesses, thus greatly reduces the computational cost of the developed maneuver detection method.

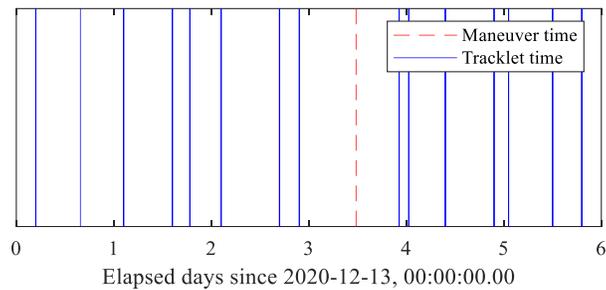

(a) Sentinel-3A

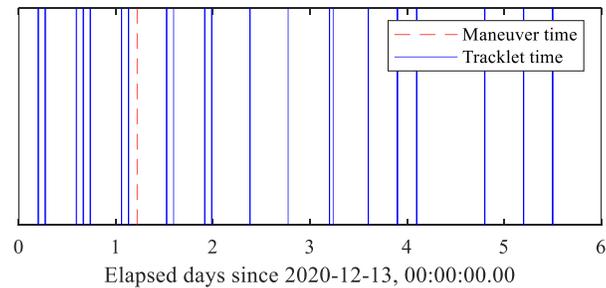

(b) Sentinel-6A

**Fig. 7.** Temporal distribution of simulated optical tracklets

*4.2. Result Analysis*

For the Sentinel-3A spacecraft, the first 5 post-maneuver tracklets arrive at $t_f+10.30\,\mathrm{h}$, $t_f+12.68\,\mathrm{h}$, $t_f+21.67\,\mathrm{h}$, $t_f+33.67\,\mathrm{h}$, and $t_f+37.23\,\mathrm{h}$. The latest update time of the pre-maneuver orbit is $t_0=t_f-14.25$ h. The first two, three, four, and five post-maneuver tracklets are used for maneuver detection, respectively. Initial orbits are estimated with post-maneuver tracklets. This initial orbit is successfully correlated to the pre-maneuver orbit with Eq. (13). The computed Mahalanobis distances with Eq. (13) are shown in Fig. 8. It is seen that the number of possible maneuver epochs ($\chi(t) \leq 3.38$) reduces when the number of post-maneuver tracklets increases.

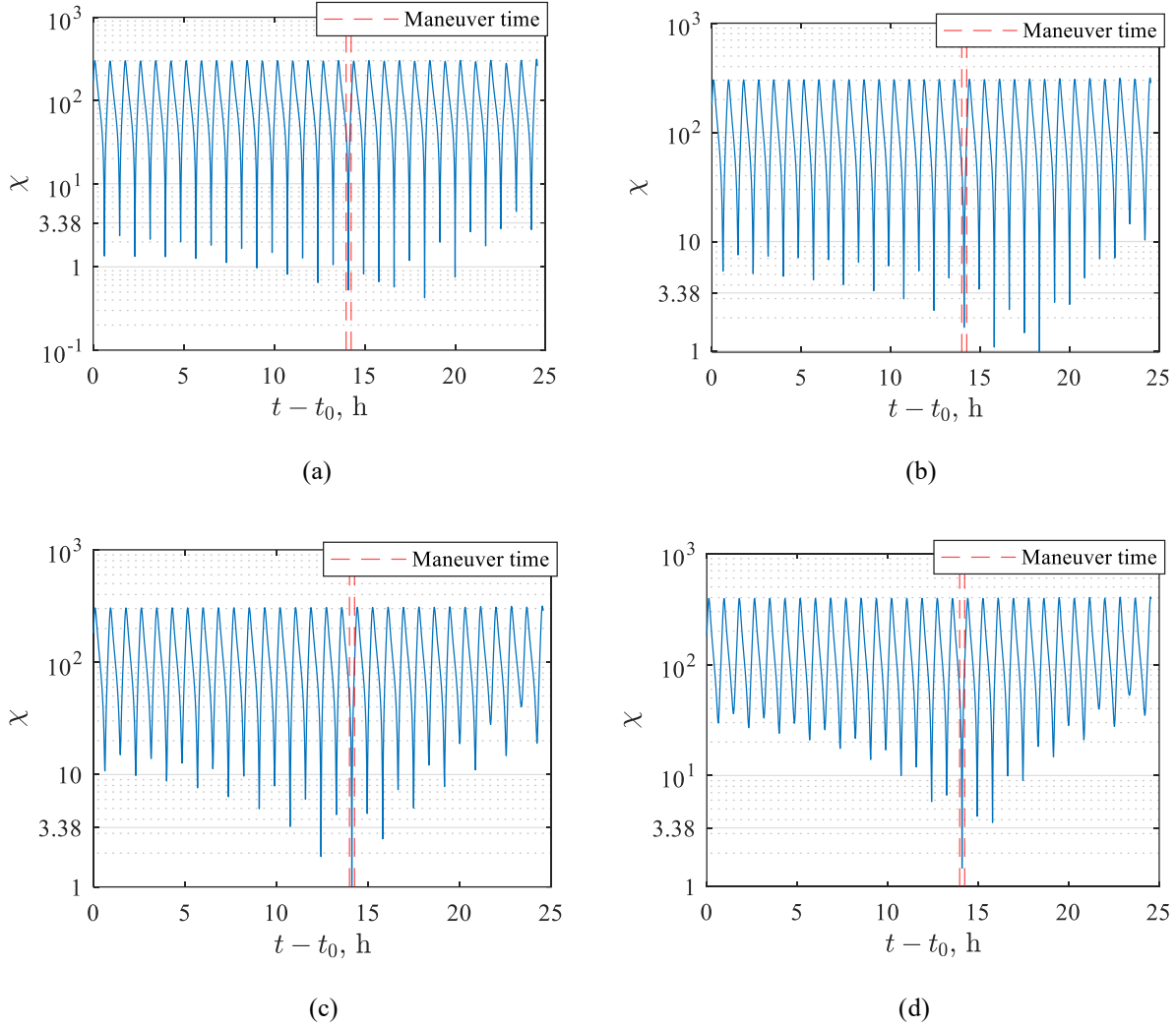

(a)

(b)

(c)

(d)

**Fig. 8.** The computed Mahalanobis distances for Sentinel-3A maneuver detection with Eq. (13) between the pre-maneuver orbit and determined post-maneuver orbit. (a) with the first two post-maneuver tracklets. (b) with the first three post-maneuver tracklets. (c) with the first four post-maneuver tracklets. (d) with the first five post-maneuver tracklets.

According to the observability analysis, the constant-thrust maneuver along the orbital normal direction has poor observability, and can also be equivalently regarded as an impulsive maneuver. The maneuver estimation results with different number of post-maneuver tracklets are shown in Fig. 9. In Fig. 9(a) ~ Fig. 9(c), multiple local minima of $J$ are found. Local minimum with the smallest $\Delta V$ is selected out as the final result based on the optimal control principle. The maneuver duration is not detectable in this case and equivalent impulsive maneuver detection results are obtained. With the first two, three, four, and five post-maneuver tracklets, the estimated impulsive maneuver

information is $\Delta V = 2.257$ m/s at {2020-12-15, 22:17:46.00}, $\Delta V = 2.251$ m/s at {2020-12-16, 08:24:34.16}, $\Delta V = 2.263$ m/s at {2020-12-16, 10:05:40.88}, and $\Delta V = 2.290$ m/s at {2020-12-16, 11:46:51.92}, respectively.

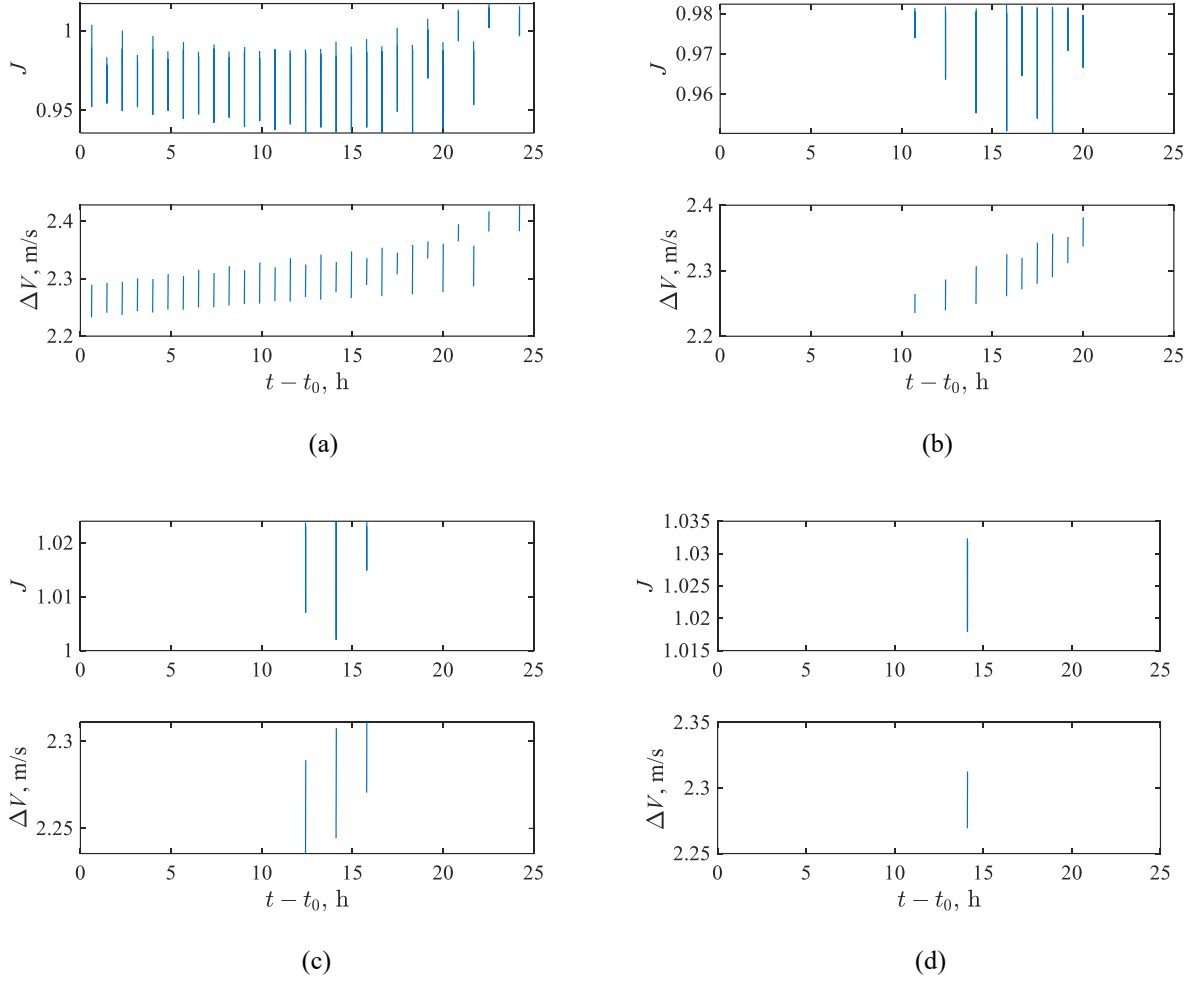

Fig. 9. Maneuver detection results for the Sentinel-3A spacecraft with different number of post-maneuver tracklets. (a) with the first two post-maneuver tracklets. (b) with the first three post-maneuver tracklets. (c) with the first four post-maneuver tracklets. (d) with the first five post-maneuver tracklets.

For the Sentinel-6A spacecraft, the first 5 post-maneuver tracklets arrive at $t_f + 7.13$ h, $t_f + 8.86$ h, $t_f + 16.57$ h, $t_f + 18.30$ h, and $t_f + 27.74$ h. The latest updated time of the pre-maneuver orbit is $t_0 = t_f - 2.30$ h. The first two, three, four, and five post-maneuver tracklets are used for maneuver detection, respectively. Initial orbits are estimated with post-maneuver tracklets. When the first two tracklets are used, the computed Mahalanobis distances with Eq. (13) are shown in Fig. 10(a). It is seen that Mahalanobis distances less than 3.38 are found. Then, impulsive maneuver detection is conducted and the estimated maneuver information is $\Delta V = 5.166$ m/s at {2020-12-14, 05:19:35.00}. In this case, the maneuver duration can hardly be detected. When the first three, four, and five tracklets are used, the effect of maneuver duration becomes evident. The computed Mahalanobis distances with Eq. (13) are all beyond 3.38. Therefore, the Mahalanobis distances are furtherly computed with Eq. (14) and an example is shown in Fig. 11.

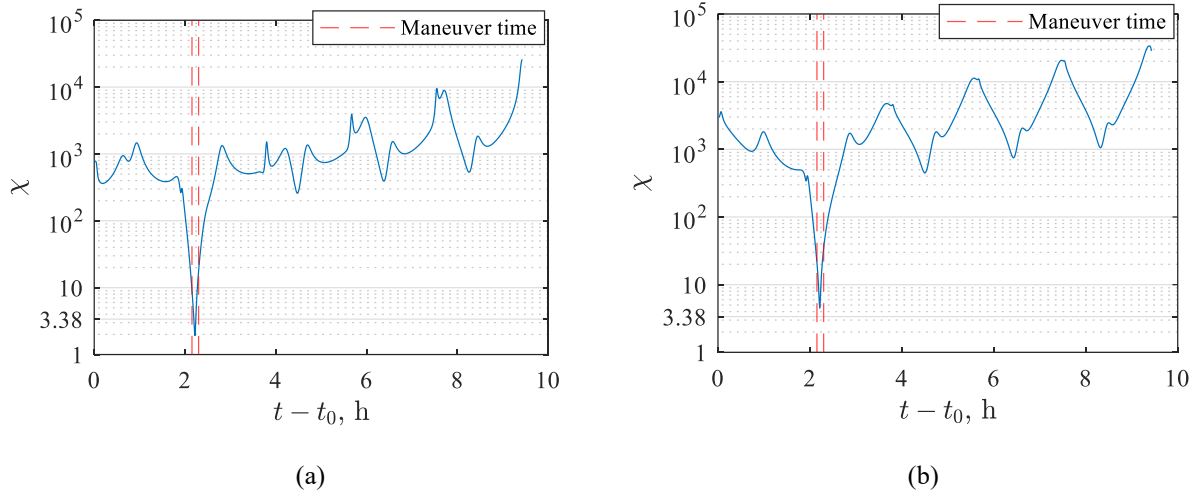

Fig. 10. The computed Mahalanobis distances for Sentinel-6A maneuver detection with Eq. (13) between the pre-maneuver orbit and determined post-maneuver orbit. (a) with the first two post-maneuver tracklets. (b) with the first three post-maneuver tracklets.

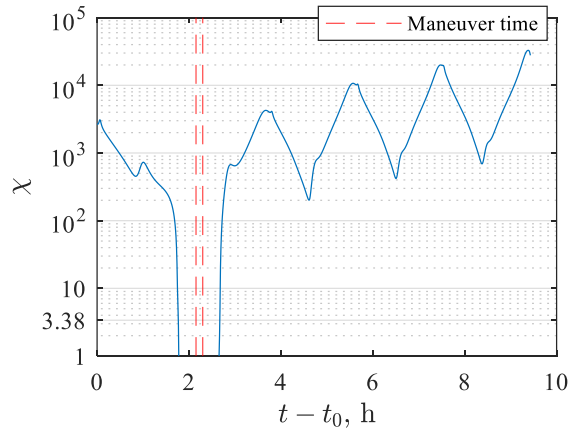

Fig. 11. The computed Mahalanobis distances for Sentinel-6A maneuver detection with Eq. (14) between the pre-maneuver orbit and determined post-maneuver orbit with the first three post-maneuver tracklets.

The long-duration maneuver detection results of the Sentinel-6A spacecraft with the first three post-maneuver tracklets is shown in Fig. 11. The estimated $\Delta V$, maneuver starting time and maneuver ending time are 5.171 m/s, {2020-12-14, 05:15:16.00}, and {2020-12-14, 05:23:16.00}, respectively. Increasing the number of post-maneuver tracklets slightly change the maneuver detection results.

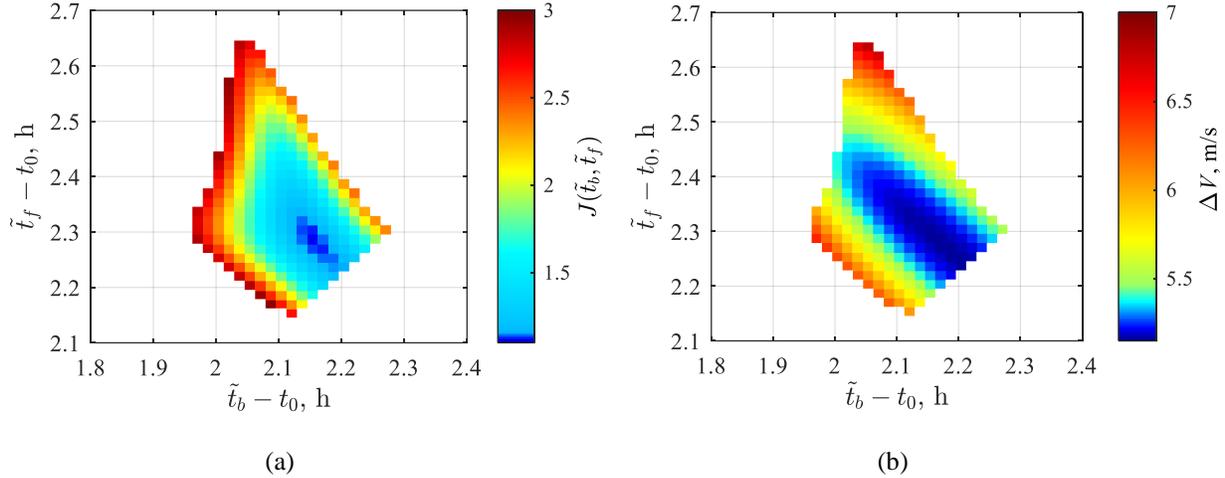

**Fig. 12.** Long-duration maneuver detection results for the Sentinel-6A spacecraft with the first three post-maneuver tracklets.

The maneuver cases listed in Table 2 are all long-duration maneuvers. But, their maneuver durations are still relatively short and can be roughly regarded as impulsive maneuvers if the accuracy requirement is low. For the Sentinel-3A spacecraft, the maneuver duration naturally has poor observability as the thrust is along the orbital normal direction. For the Sentinel-6A spacecraft, if the maneuver duration is artificially increased to 1800 s (i.e., the maneuver ending time is changed to 2020-12-14, 05:45:42 and $\Delta V = 17.716$ m/s) and other conditions stay the same, the maneuver detection with an impulsive maneuver model completely fails. The maneuver detection with a long-duration maneuver model has good results, which is shown in Fig. 13. With the first three post-maneuver tracklets, the estimated $\Delta V$, maneuver starting time and maneuver ending time are 17.730 m/s, {2020-12-14, 05:16:07.00}, and {2020-12-14, 05:46:07.00}, respectively.

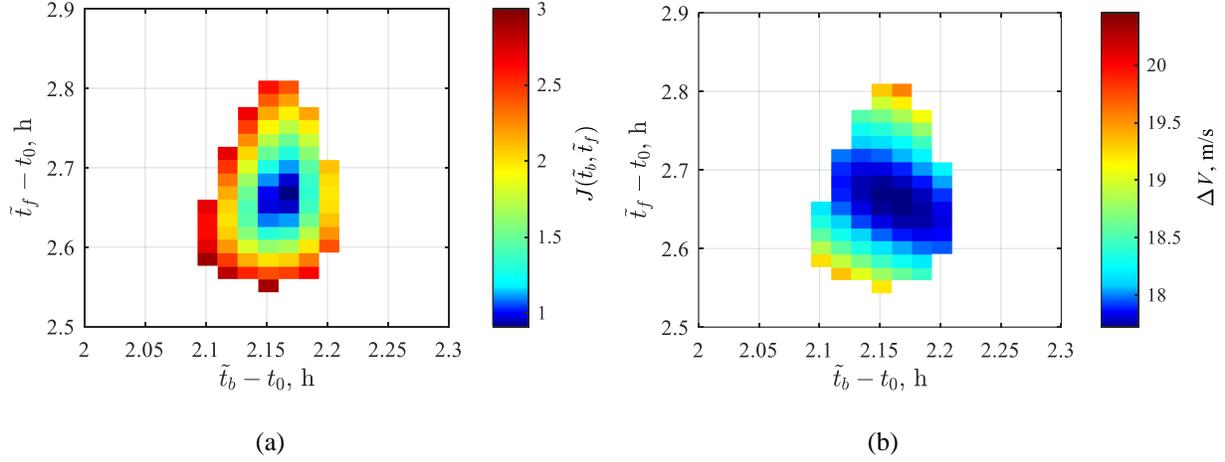

**Fig. 13.** Long-duration maneuver detection results for the Sentinel-6A spacecraft with the first three post-maneuver tracklets when the maneuver duration is increased to 1800 s.

## 5. Conclusion

In this study, a long-duration maneuver detection method in the SBSS scenario is proposed. It is shown that simplification with an impulsive maneuver model can cause significant orbital deviations when the maneuver duration is long. Existing O2O correlation method for maneuver detection is modified to deal with effects of maneuver duration. The developed maneuver detection method estimates the maneuver starting/ending times and thrust acceleration vector instead of sudden velocity change and its time. With different combinations of maneuver starting/ending time guesses, thrust acceleration vectors are estimated by nonlinear least square methods. The WRMS of angular residuals is used to judge whether corresponding maneuver policies are possible. Due to the nonlinearity of orbital dynamics there may be multiple local minima of WRMS. Acceptable local minimum with the minimum $\Delta V$ is selected as the global optimum solution according to the minimum-fuel principle. Particularly, observability analysis shows that maneuver parameters are highly coupled and poorly observable when the thrust is along the orbital normal direction. This phenomenon is analytically explained with aid of relative orbital motion. Tests with real maneuver data of the Sentinel-3A spacecraft and the Sentinel-6A spacecraft demonstrate the validity and capability of the proposed maneuver detection method. Nonetheless, the proposed method still has some limitations: 1) Only single-burn maneuver cases are considered. In data-scarce scenarios, a spacecraft may have maneuvered more than once between two subsequent optical tracklets. 2) The thrust acceleration is assumed to be constant in the orbital frame. This assumption agrees with most realistic missions, but may become invalid when the spacecraft performs delicate and sophisticated maneuvers. All in all, the developed maneuver detection method can be integrated in operational cataloging chains to increase the robustness and accuracy of uncooperative spacecraft maneuver detections.

**Appendix. The Solution Procedure for Eq. (14)**

For each epoch $t_i$, the following optimization problem is solved:

$$\chi(t_i) = \min_{\Delta r \in \mathbb{R}^3} \left(r_2(t_i) - r_1(t_i) - \Delta r\right)^T \left(P_2(t_i) + P_1(t_i)\right)^{-1} \left(r_2(t_i) - r_1(t_i) - \Delta r\right) \quad \text{s.t.} \quad \Delta r^T \Delta r \leq d^2 \quad \text{(A1)}$$

If $\left(r_2(t_i) - r_1(t_i)\right)^T \left(r_2(t_i) - r_1(t_i)\right) \leq d^2$, it is clear that $\chi(t_i) = 0$. Otherwise, compute the eigenvalue decomposition of $P_2(t_i) + P_1(t_i)$:

$$P_2(t_i) + P_1(t_i) = [\eta_1, \eta_2, \eta_3] \begin{bmatrix} \lambda_1 & 0 & 0 \\ 0 & \lambda_2 & 0 \\ 0 & 0 & \lambda_3 \end{bmatrix} [\eta_1, \eta_2, \eta_3]^T \quad \text{(A2)}$$

where $\lambda_1$, $\lambda_2$, and $\lambda_3$ are the positive eigenvalues of $P_2(t_i) + P_1(t_i)$. $\eta_1$, $\eta_2$, and $\eta_3$ are corresponding eigenvectors. Then, we have

$$r_2(t_i) - r_1(t_i) = a_1 \eta_1 + a_2 \eta_2 + a_3 \eta_3 \quad \text{(A3)}$$

$$\Delta r = b_1 \eta_1 + b_2 \eta_2 + b_3 \eta_3 \quad \text{(A4)}$$

with

$$b_1^2 + b_2^2 + b_3^2 \leq d^2 \quad \text{(A5)}$$

Furtherly, we have

$$\left(r_2(t_i) - r_1(t_i) - \Delta r\right)^T \left(P_2(t_i) + P_1(t_i)\right)^{-1} \left(r_2(t_i) - r_1(t_i) - \Delta r\right)$$
$$= \frac{(a_1 - b_1)^2}{\lambda_1} + \frac{(a_2 - b_2)^2}{\lambda_2} + \frac{(a_3 - b_3)^2}{\lambda_3} \quad \text{(A6)}$$

Equation (A1) is turned into

$$\chi(t_i) = \min_{\Delta r \in \mathbb{R}^3} \frac{(a_1 - b_1)^2}{\lambda_1} + \frac{(a_2 - b_2)^2}{\lambda_2} + \frac{(a_3 - b_3)^2}{\lambda_3} \quad \text{s.t.} \quad b_1^2 + b_2^2 + b_3^2 \leq d^2 \quad \text{(A7)}$$

Let

$$L = \frac{(a_1 - b_1)^2}{\lambda_1} + \frac{(a_2 - b_2)^2}{\lambda_2} + \frac{(a_3 - b_3)^2}{\lambda_3} + \xi(b_1^2 + b_2^2 + b_3^2 - d^2) \tag{A8}$$

where $\xi$ is the Lagrange multiplier. The KKT conditions [38] derive

$$\begin{cases} \frac{\partial L}{\partial b_1} = 0 \\ \frac{\partial L}{\partial b_2} = 0 \\ \frac{\partial L}{\partial b_3} = 0 \\ \xi(b_1^2 + b_2^2 + b_3^2 - d^2) = 0 \end{cases} \tag{A9}$$

If $\xi = 0$, it is obtained that $a_1 = b_1$, $a_2 = b_2$, and $a_3 = b_3$. This result agrees with $(r_2(t_i) - r_1(t_i))^T (r_2(t_i) - r_1(t_i)) \leq d^2$. If $\xi \neq 0$, we have $b_1^2 + b_2^2 + b_3^2 - d^2 = 0$, which derives

$$b_1 = \frac{a_1}{1 + \xi \lambda_1}; b_2 = \frac{a_2}{1 + \xi \lambda_2}; b_3 = \frac{a_3}{1 + \xi \lambda_3} \tag{A10}$$

$$\left(\frac{a_1}{1 + \xi \lambda_1}\right)^2 + \left(\frac{a_2}{1 + \xi \lambda_2}\right)^2 + \left(\frac{a_3}{1 + \xi \lambda_3}\right)^2 = d^2 \tag{A11}$$

$\xi$ can be easily solved by numerical methods according to Eq. (A11), such as the fzero function in MATLAB.

## Declaration of Competing Interest

The authors declare that they have no known competing financial interests or personal relationships that could have appeared to influence the work reported in this paper.

## Acknowledgment

This research was supported by National Key Research and Development Program of China (grant number 2020YFC1511700).

**References**

1. Liu Y, Chi R, Pang B, et al. Space debris environment engineering model 2019: Algorithms improvement and comparison with ORDEM 3.1 and MASTER-8. *Chin J Aeronaut* 2024;**37**(5):392−409. https://doi.org/10.1016/j.cja.2023.12.004

2. NASA Orbital Debris Program Office. Space missions and satellite box score. *Orbital Debris Quarterly News* 2024;**28**(4);11. Available from: https://www.orbitaldebris.jsc.nasa.gov/quarterly-news/pdfs/ODQNv28i4.pdf

3. Cai H, Xue C, Sun X, et al. Multi-sensor possibility PHD filter for space situational awareness, *Chin J Aeronaut* 2024. https://doi.org/10.1016/j.cja.2024.08.026.

4. Space Exploration Technologies Corp. SpaceX constellation status report: December 1, 2023 – May 31, 2024 [Internet]. 2024 Nov. https://licensing.fcc.gov/cgi-bin/ws.exe/prod/ib/forms/reports/related_filing.hts?f_key=-443498&f_number=SATMOD2020041700037

5. Goff GM. Orbit estimation of non-cooperative maneuvering spacecraft [dissertation]. Ohio: Air Force Institute of Technology; 2015

6. Lubey DP. Maneuver detection and reconstruction in data sparse systems with an optimal control based estimator [dissertation]. Colorado: University of Colorado; 2015

7. Pastor A, Escribano G, Sanjurjo-Rivo M, Escobar D. Satellite maneuver detection and estimation with optical survey observations. *J Astronaut Sci* 2022; **69**:879–917. https://doi.org/10.1007/s40295-022-00311-5

8. Porcelli L, Pastor A, Cano A, et al. Satellite maneuver detection and estimation with radar survey observations. *Acta Astronaut* 2022:**201**:274–287. https://doi.org/10.1016/j.actaastro.2022.08.021

9. Sharma J, Stokes GH, von Braun C, et al. Toward operational space-based space surveillance. *Linc Lab J* 2022;**13**(2):309–334

10. Silha J, Schildknecht T, Hinze A, et al. Capability of a space-based space surveillance system to detect and track objects in GEO, MEO and LEO orbits. *Proceedings of 65th International Astronautical Congress*. Toronto: IAF; 2014

11. Cheng Q, Zhang W. Scientific issues and critical technologies in planetary defense. *Chin J Aeronaut* 2024;**37**(11):24−65. https://doi.org/10.1016/j.cja.2024.07.004

12. Escribano G, Sanjurjo-Rivo M, Siminski JA, et al. Automatic maneuver detection and tracking of space objects in optical survey scenarios based on stochastic hybrid systems formulation. *Adv Space Res* 2022;**69**:3460–3477. https://doi.org/10.1016/j.asr.2022.02.034

13. Aaron BS. Geosynchronous satellite maneuver detection and orbit recovery using ground based optical tracking [dissertation]. Cambridge: Massachusetts Institute of Technology; 2006

14. Folcik ZJ, Cefola PJ, Abbot RI. Geo maneuver detection for space situational awareness. *Adv Astronaut Sci* 2008;**129**:523−550

15. Kelecy T, Jah M. Detection and orbit determination of a satellite executing low thrust maneuvers. *Acta Astronaut* 2010;**66**:798–809. https://doi.org/10.1016/j.actaastro.2009.08.029

16. Goff GM, Black JT, Beck JA. Orbit estimation of a continuously thrusting spacecraft using variable dimension filters. *J Guid Control Dynam* 2015;**38**(12):2407–2420. https://doi.org/10.2514/1.G001091

17. Goff GM, Showalter D, Black JT, Beck JA. Parameter requirements for noncooperative satellite maneuver reconstruction using adaptive filters. *J Guid Control Dynam* 2015;**38**(3):361–374. https://doi.org/10.2514/1.G000941

18. Co HC, Scheeres DJ. Maneuver detection with event representation using thrust fourier coefficients. *J Guid Control Dynam* 2006;**39**(5):1080–1091. https://doi.org/10.2514/1.G001463

19. Guang Z, Bi X, Zhao H, Liang B. Non-cooperative maneuvering spacecraft tracking via a variable structure estimator. *Aerosp Sci Technol* 2018;**79**:352–363. https://doi.org/10.1016/j.ast.2018.05.052

20. Montilla JM, Sanchez JC, Vazquez R, et al. Manoeuvre detection in Low Earth Orbit with radar data. *Adv Space Res* 2023;**72**:2689–2709. https://doi.org/10.1016/j.asr.2022.10.026

21. Conway BA. Spacecraft trajectory optimization. New York: Cambridge University Press; 2010. p. 1–35.

22. Holzinger MJ, Scheeres DJ, Alfriend KT. Object correlation, maneuver detection, and characterization using control-distance metrics. *J Guid Control Dynam* 2012;**35**(4):1312–1325. https://doi.org/10.2514/1.53245

23. Singh N, Horwood JT, Poore AB. Space object maneuver detection via a joint optimal control and multiple hypothesis tracking approach. *Proceedings of the 22nd AAS/AIAA Space Flight Mechanics Meeting*. Charleston: AIAA; 2012

24. Lubey DP, Scheeres DJ. Identifying and estimating mismodeled dynamics via optimal control policies and distance metrics. *J Guid Control Dynam* 2014;**37**(5):1512–1523. https://doi.org/10.2514/1.G000369

25. da Graça Marto S, Díaz Riofrío S, Ilioudis C, et al. Satellite Manoeuvre Detection with Multistatic Radar. *J Astronaut Sci* 2023;**70**:36. https://doi.org/10.1007/s40295-023-00399-3

26. Gao Y. Direct optimization of low-thrust many-revolution Earth-orbit transfers. *Chin J Aeronaut* 2009;**22**:426−433. https://doi.org/10.1016/S1000-9361(08)60121-1

27. Jiang R, Yang M, Wang S, et al. High-precision shape approximation low-thrust trajectory optimization method satisfying bi-objective index. *Chin J Aeronaut* 2022;**35**(1):436−457. https://doi.org/10.1016/j.cja.2020.11.022

28. Parrish N, Parker JS, Bradley BK. GEO observability from Earth-Moon libration point orbits. *Proceedings of the 23rd AAS/AIAA Space Flight Mechanics Meeting.* Kauai, Hawaii: AIAA; 2013

29. Sharma J, von Braun C, Gaposchkin EM. Space-based visible data reduction. *J Guid Control Dynam* 2000;**23**(1):170–174. https://doi.org/10.2514/2.4507

30. DeMars KJ. Nonlinear orbit uncertainty prediction and rectification for space situational awareness [dissertation]. Texas: Texas A&M University; 2010

31. Serra R, Yanez C, Frueh C. Tracklet-to-orbit association for maneuvering space objects using optimal control theory. *Acta Astronaut* 2021;**181**:271–281. https://doi.org/10.1016/j.actaastro.2021.01.026

32. Siminski JA, Montenbruck O, Fiedler H, Schildknecht T. Short-arc tracklet association for geostationary objects. *Adv Space Res* 2014:**53**(8):1184–1194. https://doi.org/10.1016/10.1016/j.asr.2014.01.017

33. Cai H, Yang Y, Gehly S, et al. Improved tracklet asociation for space objects using short-arc optical measurements. *Acta Astronaut* 2018;**151**:836–847. https://doi.org/10.1016/10.1016/j.actaastro.2018.07.024

34. Hill K, Alfriend KT, Sabol C. Covariance-based uncorrelated track association. *AIAA/AAS Astrodynamics Specialist Conference and Exhibition*. Austin: AIAA; 2008. https://doi.org/10.2514/6.2008-7211

35. Tapley BD, Schutz BE, Born GH. Statistical Orbit Determination. London: Elsevier; 2004. p. 176–194

36. Chen P, Mao X, Han J, Sun X. Observability Analysis for orbit determination using spaceborne gradiometer. *J Aerosp Eng* 2023;**36**(2):04022122. https://doi.org/10.1061/JAEEEZ.ASENG-4619

37. Vallado DA. Fundamentals of astrodynamics and applications, 4th ed. California: Microcosm Press; 2013. p. 538–584

38. Boyd S, Vandenberghe L. Convex Optimization. New York: Cambridge University Press; 2004. p. 243–244